\documentclass[12pt,a4paper]{amsart}
\usepackage{amssymb}
\usepackage[mathcal]{eucal}

\theoremstyle{plain}
\newtheorem{thm}{Theorem}[section]
\newtheorem{lem}[thm]{Lemma}
\newtheorem{pro}[thm]{Proposition}
\newtheorem{cor}[thm]{Corollary}

\theoremstyle{remark}
\newtheorem{rem}[thm]{Remark}
\newtheorem{que}[thm]{Question}

\numberwithin{equation}{section}

\newcommand{\N}{\mathbb{N}}
\newcommand{\Z}{\mathbb{Z}}
\newcommand{\Q}{\mathbb{Q}}
\newcommand{\F}{\mathbb{F}}

\newcommand{\SD}{D^\textup{semi}}

\DeclareMathOperator{\rk}{rk}
\DeclareMathOperator{\Cen}{C}
\DeclareMathOperator{\Zen}{Z}
\DeclareMathOperator{\GL}{GL}
\DeclareMathOperator{\Sym}{Sym}

\begin{document}
\title[On the rank of compact $p$-adic Lie groups]{On the rank of compact
  $p$-adic Lie groups}


\author{Benjamin Klopsch} \address{Department of Mathematics, Royal
  Holloway, University of London, Egham TW20 0EX, United Kingdom}
\email{Benjamin.Klopsch@rhul.ac.uk}

\begin{abstract}
  The rank of a profinite group $G$ is the basic invariant $\rk(G) :=
  \sup \{ d(H) \mid H \leq G \}$, where $H$ ranges over all closed
  subgroups of $G$ and $d(H)$ denotes the minimal cardinality of a
  topological generating set for $H$.  A compact topological group $G$
  admits the structure of a $p$-adic Lie group if and only if it
  contains an open pro-$p$ subgroup of finite rank.  

  For every compact $p$-adic Lie group $G$ one has $\rk(G) \geq
  \dim(G)$, where $\dim(G)$ denotes the dimension of $G$ as a $p$-adic
  manifold.  In this paper we consider the converse problem, bounding
  $\rk(G)$ in terms of $\dim(G)$.

  Every profinite group $G$ of finite rank admits a maximal finite
  normal subgroup, its periodic radical $\pi(G)$.  One of our main
  results is the following.  Let $G$ be a compact $p$-adic Lie group
  such that $\pi(G) = 1$, and suppose that $p$ is odd.  If $\{ g \in G
  \mid g^{p-1} = 1 \}$ is equal to $\{1\}$, then $\rk(G) = \dim(G)$.
\end{abstract}

\keywords{compact $p$-adic Lie group, rank, dimension, number of
  generators}

\subjclass[2000]{20E18, 22E20}

\maketitle

\section{Introduction}

In the $1960$s Lazard developed a sophisticated theory of $p$-adic Lie
groups.  One of the high points of his work was the algebraic
characterisation of such groups, providing a solution to the $p$-adic
analogue of Hilbert's $5$th problem.  A more recent approach to the
theory of compact $p$-adic Lie groups, inspired by Lazard's work and
focusing on group-theoretic aspects, is given in~\cite{DiDuMaSe99}.

The \emph{rank} of a profinite group $G$ is the basic invariant
$\rk(G) := \sup \{ d(H) \mid H \leq G \}$, where $H$ ranges over all
closed subgroups of $G$ and $d(H)$ denotes the minimal cardinality of
a topological generating set for $H$.  According to
\cite[Corollary~8.34]{DiDuMaSe99}, a topological group $G$ admits the
structure of a compact $p$-adic Lie group if and only if it is a
profinite group of finite rank which contains an open pro-$p$
subgroup.  Moreover, the analytic structure which makes $G$ into a
$p$-adic Lie group, whenever it exists, is unique.  See
\cite[Interlude~A]{DiDuMaSe99} for a wide range of alternative
characterisations of pro-$p$ groups of finite rank.

Another key invariant of a compact $p$-adic Lie group $G$ is its
\emph{dimension}, $\dim(G)$, as a $p$-adic manifold.  Algebraically,
$\dim(G)$ can be described as $d(U)$ where $U$ is any uniformly
powerful open pro-$p$ subgroup of $G$.  Consequently, for every
compact $p$-adic Lie group $G$ one has $\rk(G) \geq \dim(G)$.
Uniformly powerful pro-$p$ groups are rather special $p$-adic analytic
groups which are torsion-free and in many ways behave like lattices.
For instance, one has $\rk(U) = \dim(U)$ for every uniformly powerful
pro-$p$ group $U$.  

Straightforward examples show that, in general, to bound the rank
$\rk(G)$ of a compact $p$-adic Lie group $G$ in terms of $\dim(G)$,
one needs to restrain the torsion of $G$.  In the first instance, it
is thus natural to consider groups which are torsion-free.

\begin{que} \label{que:torsion-free} Is it true that for every
  torsion-free compact $p$-adic Lie group $G$ one has $\rk(G) =
  \dim(G)$?
\end{que}

In Section~\ref{sec:Laffey} we show that a result of
Laffey~\cite{La73}, which bounds the number of generators of finite
$p$-groups, readily implies
  
\begin{pro}\label{pro:Laffey}
  Let $G$ be a torsion-free compact $p$-adic Lie group, and suppose
  that $p$ is odd.  Then $\rk(G) = \dim(G)$.
\end{pro}

It is remarkable that a $p$-group argument yields such a result, but
unfortunately the method of proof does not appear to lead further.
Proposition~\ref{pro:Laffey} is restricted to groups without any
torsion and falls short of dealing with $2$-adic Lie groups.  An
analogue of Laffey's theorem for $p=2$ is given in
\cite[Corollary~3.10]{GoKl10}; it has, however, rather weak
implications in the present context.

\medskip

In the present paper we take a new approach which leads to the more
general Theorems~\ref{thm:A} and~\ref{thm:B} below and which, we
think, will be more suitable for dealing with $2$-adic Lie groups.
Indeed, our results are formulated for compact $p$-adic Lie groups
which possess no non-trivial finite normal subgroups, a condition
which also arises in the study of lattices.  Every profinite group $G$
of finite rank admits a maximal finite normal subgroup $\pi(G)$, the
\emph{periodic radical} of~$G$.  If $G$ is $p$-adic analytic, so is
$G/\pi(G)$ and $\dim(G) = \dim(G/\pi(G))$.  These facts provide a
natural reduction to groups $G$ satisfying $\pi(G) = 1$.  For every
prime~$\ell$, we define the \emph{$\ell$-rank} of a profinite group
$G$ to be $\rk_\ell(G) := \rk(S)$ where $S$ is any Sylow pro-$\ell$
subgroup of~$G$.  Our main result is

\begin{thm}\label{thm:A}
  Let $G$ be a compact $p$-adic Lie group with $\pi(G) = 1$, and
  suppose that $p$ is odd.  Then one has
  $$
  \max \{ \rk_\ell(G) \mid \text{$\ell>2$ prime} \} = \rk_p(G) =
  \dim(G)
  $$
  and
  $$
  \rk_2(G) \leq
  \begin{cases}
    \lfloor 3 \dim(G) / 2 \rfloor & \text{if $p \equiv_4 1$,} \\
    \dim(G) & \text{if $p \equiv_4 3$.}
  \end{cases}
  $$
\end{thm}

A theorem of Guralnick and Lucchini, generalising work of Kov\'acs,
shows that every profinite group $G$ satisfies $\rk(G) \leq \sup \{
\rk_\ell(G) \mid \text{$\ell$ prime} \} + 1$; see~\cite{Gu89,Lu89}.
This and related results allow us to prove
 
\begin{thm}\label{thm:B}
  Let $G$ be a compact $p$-adic Lie group with $\pi(G) = 1$, and
  suppose that $p$ is odd.  Then one has
  \begin{enumerate}
  \item[(1)] $\rk(G) \leq \max \{ \dim(G), \rk_2(G) \} + 1$;
  \item[(2)] $\rk(G) \leq \max \{ \dim(G)+1, \rk_2(G) \}$ if $G$ is
    prosoluble;
  \item[(3)] $\rk(G) \leq \dim(G)+1$ if $p \equiv_4 3$;
  \item[(4)] $\rk(G) = \dim(G)$ if $G$ has trivial $(p-1)$-torsion, i.e.\
    if $\{ g \in G \mid g^{p-1} = 1 \}$ is equal to $\{1\}$.
  \end{enumerate}
\end{thm}

\begin{rem}\label{rem:beispiele}
  Relatively simple examples show that Theorems~\ref{thm:A} and
  \ref{thm:B} are to some extent best possible.  Indeed, for $p
  \equiv_4 1$ the matrix group
  $$
  S := \langle
  \begin{pmatrix}
    0 & 1 \\ 1 & 0
  \end{pmatrix},
  \begin{pmatrix}
    \sqrt{-1} & 0 \\ 0 & \sqrt{-1}
  \end{pmatrix},
  \begin{pmatrix}
    -1 & 0 \\ 0 & 1
  \end{pmatrix}
  \rangle \leq \GL_2(\Z_p)
  $$
  has order $16$, acts irreducibly on the abelian lattice $\Z_p^2$ and
  requires a minimum of $3$ generators.  Thus for any $r \in \N$ and any
  prime $p \equiv_4 1$, the $2r$-dimensional compact $p$-adic Lie
  group $G = (S \ltimes \Z_p^2)^r$ has $\pi(G) = 1$ and $\rk_2(G) = 3r
  = 3 \dim(G)/2$.

  For any prime $p$, any divisor $m \neq 1$ of $p-1$ and any $d \in
  \N$, the group $G = T \ltimes \Z_p^d$, where $T \cong C_m$ consists
  of the scalar matrices in $\GL_d(\Z_p)$ corresponding to $m$th roots
  of unity in $\Z_p$, has $\pi(G) = 1$ and $\rk(G) = d+1 = \dim(G)
  +1$.
\end{rem}

While the focus of the present paper is on $p$-adic Lie groups for odd
primes~$p$, it remains a challenging open problem to describe the
precise relation between rank and dimension for compact $2$-adic Lie
groups.  Basic examples indicate that things are somewhat different
for the prime $2$.  For instance, for every $r \in \N$ the $r$-fold
direct power of the pro-$2$ dihedral group $D = C_2 \ltimes \Z_2$
satisfies $\pi(D^r) = 1$ and $\rk(D^r) = 2r = 2 \dim(D^r)$.

We recall that a pro-$p$ group $G$ is said to be \emph{$d$-maximal} if
$d(H) < d(G)$ for every proper open subgroup $H$ of $G$.  It is known
that, for odd $p$, every $d$-maximal finite $p$-group is nilpotent of
class at most~$2$; cf.~\cite{GoKl10} and references therein.  Our
proof of Theorem~\ref{thm:A} indicates that
Question~\ref{que:torsion-free} for $p=2$ is linked to the long
standing problem of understanding the structure of $d$-maximal finite
$2$-groups and pro-$2$ groups.  Of particular relevance would be a
positive answer to

\begin{que}\label{que:d-maximal}
  Is it true that every $d$-maximal pro-$2$ group is soluble?
\end{que}

Besides using the concept of $d$-maximal groups, we employ in our
proof of Theorem~\ref{thm:A} standard techniques from the theory of
$p$-adic analytic pro-$p$ groups.  Furthermore we take advantage of
the description of maximal $\ell$-subgroups of general linear groups
over $\Q_p$, given in \cite{LePl86}, and the classification of
indecomposable $\Z_p C$-modules for cyclic groups $C$ of order $p$,
provided in \cite{HeRe62}.  Of independent interest is our elementary
proof of the following auxiliary result, which one may compare with
more general, asymptotic estimates for the minimal number of
generators of finite linear groups; for instance, see
\cite{Is72,DiKo91,KoRo91,LuMeMo01}.

\begin{thm}\label{thm:GL_d_F_p}
  Suppose that $p$ is odd.  Let $d \in \N$, and let $\ell$ be prime
  with $\ell \neq p$.  Let $d_0 \in \{0,1\}$ such that $d \equiv_2
  d_0$, and set $m(p,\ell) := \min \{ n \in \N \mid \ell \text{
    divides } p^n - 1 \}$.  Then
  $$
  \rk_\ell(\GL_d(\Z_p)) = \rk_\ell(\GL_d(\F_p)) = 
  \begin{cases}
    \lfloor d / m(p,\ell) \rfloor & \text{if $p \not \equiv_\ell 1$,} \\
    (3d - d_0)/2 & \text{if $\ell = 2$ and $p
      \equiv_4 1$,} \\
    d & \text{in all other cases.}
  \end{cases}
  $$
\end{thm}

\smallskip

\noindent \textit{Notation.}  Throughout, $p$ and $\ell$ denote
primes.  The field of $p$-adic numbers and the ring of $p$-adic
integers are denoted by $\Q_p$ and $\Z_p$.  The commutator subgroup of
a group $G$ is denoted by $[G,G]$.  The centre of $G$ is denoted by
$\Zen(G)$, the centraliser in $G$ of a subset $S \subseteq G$ by
$\Cen_G(S)$.  The minimal cardinality of a generating set for a group
$G$ is denoted by $d(G)$.  Likewise the minimal cardinality of a
generating set for a module $M$ over a ring $R$ is denoted by
$d_R(M)$.  For $n \in \N_0$ it is customary to denote by $G^n$ the
$n$-fold direct power of $G$; in a few places, we use the same
notation, $G^n$ for $n$ a power of $p$, to denote the characteristic
subgroup of $G$ generated by all $n$th powers of elements of~$G$.  All
wreath products considered are permutational wreath products.

When $G$ is a topological group, subgroups $H$ of $G$ are tacitly
taken to be closed, $d(G)$ tacitly refers to the minimal number of
topological generators, etc.  A key concept in the theory of $p$-adic
Lie groups is that of a powerful pro-$p$ group.  We refer to
\cite[Part~I]{DiDuMaSe99} for standard properties of uniformly
powerful pro-$p$ groups.  A profinite group $G$ is said to be
\emph{just infinite} if it is infinite and if every proper quotient of
$G$ is finite.


\section{A short proof using Laffey's bound}
\label{sec:Laffey}

Clearly, Proposition~\ref{pro:Laffey} is a consequence of
Theorem~\ref{thm:A}.  But it also admits a short independent proof,
based on a result of Laffey~\cite{La73}.

\begin{proof}[Proof of Proposition~\ref{pro:Laffey}] For every prime
  $\ell$ with $\ell \neq p$, the Sylow pro\nobreakdash-$\ell$ subgroup
  of $G$ is finite.  Since $G$ is torsion-free, we conclude that $G$
  is a pro-$p$ group.  Let $U$ be a uniformly powerful open normal
  subgroup of $G$.  Then $d(U) = \dim(G)$ and descending to an open
  subgroup, if necessary, it suffices to show that $d(G) \leq
  \dim(G)$.  The open subgroups $U^{p^n}$, $n \in \N$, provide a base
  for the neighbourhoods of $1$ in~$G$.  Consequently, since $G$ is
  torsion-free, we find $m \in \N$ such that $\{ x \in G \mid x^p \in
  U^{p^m} \} \subseteq U$.  (In fact, a short argument shows that
  already $m=1$ works.)  As $U$ is uniformly powerful, this implies
  that for all $n \in \N$ with $n \geq m$,
  $$
  \Omega_1(G/U^{p^n}) = \langle x U^{p^n} \in G/U^{p^n} \mid x^p \in
  U^{p^n} \rangle = U^{p^{n-1}} / U^{p^n}
  $$
  has size $p^{\dim(U)}$.  Since $p$ is odd, we deduce from
  \cite[Corollary~2]{La73} that
  $$
  d(G) = \limsup_{n \in \N} d(G/U^{p^n}) \leq \limsup_{n \in \N}
  \log_p \lvert \Omega_1(G/U^{p^n}) \rvert = \dim(G).
  $$
\end{proof}


\section{Maximal $\ell$-subgroups of $\GL_d(\Q_p)$}

For $r \in \N_0$ we denote by $W_r(\ell)$ the $r$-fold iterated
permutational wreath product $C_\ell \wr \cdots \wr C_\ell$ of cyclic
groups of order $\ell$.  This finite permutation group of degree
$\ell^r$ has order $\lvert W_r(\ell) \rvert = \ell^{(\ell^r - 1)/(\ell
  - 1)}$ and requires $d(W_r(\ell)) = r$ generators.  Recall that the
Sylow $\ell$-subgroups of a symmetric group $\Sym(n)$ of finite degree
$n$ can be described in terms of iterated wreath products: each of
them is isomorphic to $\prod_{i=1}^t W_i(\ell)^{n_i}$, where $n =
\sum_{i=0}^t n_i \ell^i$ is the $\ell$-adic expansion of $n$ with
coefficients $0 \leq n_i < \ell$.

The Sylow $\ell$-subgroups of finite symmetric groups feature in the
classification of maximal $\ell$-subgroups of general linear groups;
see~\cite{LePl86}.  We define the following numerical invariants:
\begin{align*}
  m(p,\ell) & :=
  \begin{cases}
    \min \{ n \in \N \mid \ell \text{ divides } p^n - 1 \} & \text{if
      $p \neq \ell$,} \\
    p-1 & \text{if $p=\ell$,}
  \end{cases}
  \\
  a(p,\ell) & := 
  \begin{cases}
    \max \{ b \in \N_0 \mid \ell^b \text{ divides } p^{m(p,\ell)} - 1
    \} & \text{if $p \neq \ell$,} \\
    1 & \text{if $p=\ell$.}
  \end{cases}
\end{align*}
In addition, we define for $p \equiv_4 3$ (and $\ell = 2$),
$$ 
c(p,2) := \max \{ b \in \N \mid 2^b \text{ divides } p^2 - 1 \}.
$$
Furthermore, we recall that for $c \in \N$ with $c \geq 3$ the
semidihedral group of order $2^{c+1}$ is given by the presentation
$$
\SD_{2^{c+1}} = \langle x,y \mid x^2 = y^{2^c} = 1, \, y^x =
y^{-(1+2^{c-1})} \rangle.
$$

We require the following information about maximal $\ell$-subgroups of
general linear groups over the $p$-adic field $\Q_p$.

\begin{pro}\label{pro:max_p_subs}
  Let $d \in \N$ and suppose that $(p,\ell) \not = (2,2)$.  Then the
  group $\GL_d(\Q_p)$ has exactly one conjugacy class of maximal
  $\ell$-subgroups.  Let $L$ be one of the maximal $\ell$-subgroups of
  $\GL_d(\Q_p)$.

  \textup{(1)} If $\ell$ is odd, then $L \cong C_{\ell^a} \wr S$,
  where $a = a(p,\ell)$ and $S$ is a Sylow $\ell$-subgroup of
  $\Sym(\lfloor d / m(p,\ell) \rfloor)$.

  \textup{(2)} Suppose that $\ell = 2$ and that $p$ is odd.  If $p
  \equiv_4 1$, then $L \cong C_{2^a} \wr S$, where $a = a(p,2)$ and
  $S$ is a Sylow $2$-subgroup of $\Sym(d)$.  If $p \equiv_4 3$ and $d
  \equiv_2 0$, then $L \cong \SD_{2^{c+1}} \wr S$, where $c = c(p,2)$
  and $S$ is a Sylow $2$-subgroup of $\Sym(d/2)$.  If $p \equiv_4 3$
  and $d \equiv_2 1$, then $L \cong C_2 \times (\SD_{2^{c+1}} \wr S)$,
  where $c = c(p,2)$ and $S$ is a Sylow $2$-subgroup of
  $\Sym((d-1)/2)$.
\end{pro}

\begin{proof}
  The proposition is a special instance of Theorems~II.4, IV.2 and
  IV.3 in~\cite{LePl86}.  Indeed, $m(p,\ell)$ is equal to the degree
  of the $\ell$th cyclotomic field $\Q_p(\zeta_\ell)$ over $\Q_p$,
  where $\zeta_\ell$ denotes a primitive $\ell$th root of unity, and
  $a(p,\ell) = \max \{ b \in \N \mid \Q_p(\zeta_\ell) \text{ contains
    a primitive $\ell^b$th root of $1$} \}$.  Moreover, if $\ell = 2$
  and $p \equiv_4 3$, then $c = c(p,2)$ takes the same values as the
  invariant $\gamma(\Q_p)$ used in \cite[Section~IV]{LePl86}.  We
  observe that Theorem~IV.4 in op.\ cit.\ is not used, as we do not
  discuss the case $p = \ell = 2$.
\end{proof}

Our aim is to prove Theorem~\ref{thm:GL_d_F_p}.  In view of
Proposition~\ref{pro:max_p_subs}, we consider for $r \in \N_0$ the
groups
\begin{align*}
  X_{a,r}(\ell) & := C_{\ell^a} \wr W_r(\ell), && a \in \N, \\
  Y_{c,r} := Y_{c,r}(2) & := \SD_{2^{c+1}} \wr W_r(2), && c \in \N
  \text{ with } c \geq 3.
\end{align*}

We require the following lemma.

\begin{lem}\label{lem:no_fixed_points}
  Let $n \in \N$, and suppose that $\sigma \in \Sym(n)$ is a
  permutation of $\ell$-power order acting without fixed points.  Let
  $h \in \GL_n(\Z_\ell)$ be a monomial matrix with permutation part
  corresponding to $\sigma$.  Let $M$ be a $\Z_\ell \langle h
  \rangle$-submodule of the natural $\Z_\ell \langle h \rangle$-module
  $V = \Z_\ell^n$.  Then $d_{\Z_\ell \langle h \rangle}(M) \leq 2n /
  \ell$.
\end{lem}

\begin{proof}
  We observe that $n$ is a multiple of $\ell$, and we argue by
  induction on $n$.  First suppose that $n = \ell$ so that $\sigma$ is
  an $\ell$-cycle.  We are to show that $d_{\Z_\ell \langle h
    \rangle}(M) \leq 2$.  Clearly, this is true if $\ell = 2$.  Now
  suppose that $\ell > 2$.  Then we write $\det(h) = \lambda^\ell
  (1+\mu)$, where $\lambda \in \Z_\ell^*$ and either $\mu = 0$ or $\mu
  \in \ell \Z_\ell \setminus \ell^2 \Z_\ell$.  After a change of
  basis, we may assume that the action of $h$ on the standard
  $\Z_\ell$-basis $e_1, \ldots, e_\ell$ of $V$ is given by
  $$
  e_i^h = \lambda e_{i+1} \text{ for $i \in \{1, \ldots, \ell-1\}$}
  \quad \text{and} \quad e_\ell^h = \lambda (1+\mu) e_1.
  $$
  For the purpose of bounding $d_{\Z_\ell \langle h \rangle}(M)$ we
  may as well assume that $\lambda = 1$, and we treat the two
  possibilities for $\mu$ separately.

  First consider the case $\mu = 0$.  Then $h$ is the permutation
  matrix associated to the $\ell$-cycle $\sigma$ and $V$ the
  corresponding integral permutation module.  The $\Z_\ell \langle h
  \rangle$-module $M$ is free as a $\Z_\ell$-module and decomposes as
  a direct sum of indecomposable $\Z_\ell \langle h \rangle$-modules.
  Each indecomposable summand $I$ of $M$ satisfies $d_{\Z_\ell \langle
    h \rangle}(I) = 1$ and $\dim_{\Z_\ell}(I) \in \{1,\ell-1,\ell\}$;
  see~\cite[Theorem~2.6]{HeRe62}.  Since $\dim_{\Z_\ell}(M) \leq n =
  \ell$ and since $\langle h \rangle$ does not act trivially on $M$,
  it follows that $M$ is the sum of at most two indecomposable
  $\Z_\ell \langle h \rangle$-modules, and hence $d_{\Z_\ell \langle h
    \rangle}(M) \leq 2$, as wanted.

  Next consider the case $\mu \in \ell \Z_\ell \setminus \ell^2
  \Z_\ell$.  Then the minimum polynomial of $h-1$ over $\Q_\ell$ is
  $(X+1)^\ell - (1+\mu)$, an Eisenstein polynomial.  Thus the ring
  $\Z_\ell \langle h \rangle$ is isomorphic to the ring of integers
  $\mathcal{O}$ of a totally ramified extension $K$ of $\Q_\ell$ of
  degree $\ell$.  Furthermore, the $\Z_\ell \langle h \rangle$-module
  $V$ corresponds to the $\mathcal{O}$-module $\mathcal{O}$, and $M$
  to an ideal of the principal ideal domain $\mathcal{O}$.  Thus
  $d_{\Z_\ell \langle h \rangle}(M) \leq 1$, finishing the argument
  for $n = \ell$.

  Now suppose that $n > \ell$ and that the permutation group $\langle
  \sigma \rangle$ is intransitive: $\langle \sigma \rangle \leq
  \Sym(n_1) \times \Sym(n_2)$, where $n = n_1 + n_2$ with positive
  summands.  Write $\sigma_1$ and $\sigma_2$ for the images of
  $\sigma$ under projection into $\Sym(n_1)$ and $\Sym(n_2)$
  respectively, and note that each of these permutations acts without
  fixed points.  The monomial matrix $h$ admits a corresponding block
  decomposition $h = h_1 \oplus h_2$, with the permutation part of the
  monomial matrix $h_i$ corresponding to $\sigma_i$ for $i \in
  \{1,2\}$.  Clearly, the module $V$ admits a $\Z_\ell \langle h
  \rangle$-submodule $W$ such that $V_1 := V/W \cong \Z_\ell^{n_1}$
  and $V_2 := W \cong \Z_\ell^{n_2}$ are naturally modules for
  $\Z_\ell \langle h_1 \rangle$ and $\Z_\ell \langle h_2 \rangle$.
  Setting $M_1 := (M+W)/W$ and $M_2 := M \cap W$, we deduce by
  induction that
  $$
  d_{\Z_\ell \langle h \rangle}(M) \leq d_{\Z_\ell \langle h_1 \rangle
  }(M_1) + d_{\Z_\ell \langle h_2 \rangle}(M_2) \leq 2 (n_1 + n_2) /
  \ell = 2 n / \ell.
  $$

  Finally, suppose that $n > \ell$ and that the permutation group
  $\langle \sigma \rangle$ is transitive.  Then $n$ is equal to the
  order of $\sigma$, hence an $\ell$-power: $n = \ell^k$ with $k \geq
  2$.  Clearly, $\sigma^\ell \in \Sym(n)$ acts with $\ell$ orbits of
  size $\ell^{k-1}$ and has no fixed points.  Applying the argument
  given before, we deduce that $ d_{\Z_\ell \langle h \rangle}(M) \leq
  d_{\Z_\ell \langle h^\ell \rangle}(M) \leq 2n/\ell $.
\end{proof}

\begin{pro} \label{pro:Xgroups} 
  Let $a \in \N$ and $r \in \N_0$.  Then
  $$
  \rk(X_{a,r}(\ell)) =
  \begin{cases}
    3 \cdot 2^{r-1} & \text{if $\ell = 2$, $a \geq 2$ and $r \geq 1$,} \\
    \ell^r & \text{in all other cases.}
  \end{cases}
  $$
\end{pro}

\begin{proof}
  We write $X := X_{a,r}(\ell)$ and argue by induction on $r$.  If $r
  = 0$, the group $X \cong C_{\ell^a}$ is cyclic of order $\ell^a$ and
  the assertion holds.  Now suppose that $r \geq 1$.

  First we derive a lower bound for $\rk(X)$.  Let $E :=
  E_{a,r}(\ell)$ denote the homocyclic group of exponent $\ell^a$
  which arises as the base subgroup of $X = C_{\ell^a} \wr W_r(\ell)$.
  Clearly, $\rk(X) \geq \rk(E) = \ell^r$.  Next we derive for $\ell =
  2$, $a \geq 2$ and $r \geq 1$ the stronger lower bound $\rk(X) \geq
  3 \cdot 2^{r-1}$.  Indeed, the $2$-group $S$ displayed in
  Remark~\ref{rem:beispiele} has $d(S) = 3$ and is easily embedded
  into $C_4 \wr C_2$.  Since $a \geq 2$, the group $C_4 \wr C_2$
  embeds into $C_{2^a} \wr C_2$, and the decomposition $X \cong
  (C_{2^a} \wr C_2) \wr W_{r-1}(2)$ shows that the $2^{r-1}$-fold
  direct power of $S$ embeds into $X$.  We conclude that $\rk(X) \geq
  d(S^{2^{r-1}}) = 3 \cdot 2^{r-1}$.

  It remains to give an upper bound for $\rk(X)$.  Consider an
  arbitrary subgroup $H \leq X$, and decompose $X \cong
  X_{a,r-1}(\ell) \wr C_\ell$ as $X = B \rtimes C$, where $B = B_1
  \times \ldots \times B_\ell$ with $B_i \cong X_{a,r-1}(\ell)$ for $i
  \in \{1,\ldots,\ell\}$ and $C \cong C_\ell$.  In addition, we
  decompose the homocyclic base group $E$ of $X$ into $E = E_1 \times
  \ldots \times E_\ell$, where $E_i = E \cap B_i \cong
  (C_{\ell^a})^{\ell^{r-1}}$ is the homocyclic base group of the
  wreath product $B_i \cong C_{\ell^a} \wr W_{r-1}(\ell)$ for $i \in
  \{1,\ldots,\ell\}$.

  First consider the case $H \leq B$.  Induction on $r$ yields
  \begin{equation}\label{equ:d(H)}
    \begin{split}
      d(H) \leq \sum_{i=1}^\ell \rk(B_i) & = \ell \cdot \rk(X_{a,r-1})
      \\
      & =
      \begin{cases}
        3 \cdot 2^{r-1} & \text{if $\ell = 2$, $a \geq 2$ and $r \geq
          2$,}
        \\
        \ell^r & \text{in all other cases.}
      \end{cases}
    \end{split}
  \end{equation}
  Indeed, by projecting $H$ into the first factor $B_1$ of $B$, then
  projecting the kernel $N_1$ of the first projection into the second
  factor $B_2$ of $B$, and so on, we find a descending subnormal
  series $H = N_0 \geq N_1 \geq \ldots \geq N_\ell = 1$ such that, for
  $i \in \{1,\ldots,\ell\}$, the quotient $N_{i-1}/N_i$ embeds into
  $B_i$.  Clearly, $d(H) \leq \sum_{i=1}^\ell d(N_{i-1}/N_i) \leq
  \sum_{i=1}^\ell \rk(B_i)$ and \eqref{equ:d(H)} is justified.

  It remains to consider the case where $H$ is not contained in $B$,
  and it will be enough to prove that under this hypothesis $d(H) \leq
  3 \ell^{r-1}$.  Decompose $X \cong C_{\ell^a} \wr W_r(\ell)$ as $X =
  E \rtimes W$, where $E = E_1 \times \ldots \times E_\ell$ with $E_i
  = E \cap B_i$ for $i \in \{1,\ldots,\ell\}$ as before and $W \cong
  W_r(\ell) = X_{1,r-1}(\ell)$.  By induction on $r$, we have $d(HE/E)
  \leq \rk(W) = \ell^{r-1}$.  Furthermore, if $\ell = 2$, we have $d(H
  \cap E) \leq \rk(E) = \ell^r = 2 \cdot \ell^{r-1}$, thus $d(H) \leq
  d(HE/E) + d(H \cap E) \leq 3 \ell^{r-1}$, as wanted.

  Now suppose that $\ell \geq 3$.  Choose $h \in H \setminus B$ and
  regard $H \cap E$ as a $\Z_\ell \langle h \rangle$-module.  In view
  of the inequality $d(H) \leq d(HE/E) + d_{\Z_\ell \langle h
    \rangle}(H \cap E)$, it suffices to prove that $d_{\Z_\ell \langle
    h \rangle}(H \cap E) \leq 2 \ell^{r-1}$.  Note that $h$ permutes
  cyclically the $\ell$ factors $E_i$ of $E$.  Working with a
  pre-image $M$ of $H \cap E$ in the base group of the $\ell$-adic
  group $\Z_\ell \wr W$, we apply Lemma~\ref{lem:no_fixed_points} to
  deduce that $d_{\Z_\ell \langle h \rangle}(H \cap E) \leq 2
  \ell^r/\ell = 2 \ell^{r-1}$.
\end{proof}

We remark that, for $\ell \geq 5$, one can easily extend the argument
given in the proof of Proposition~\ref{pro:Xgroups} to show that every
subgroup $H$ of $X := X_{a,r}(\ell) = C_{\ell^a} \wr W_r(\ell)$ with
$d(H) = \ell^r = \rk(X)$ is contained in the homocyclic base group $E
\cong C_{\ell^a}^{\ell^r}$.  It remains a curious open problem to
classify subgroups witnessing the rank of $X_{a,r}(\ell)$ for the
small primes $\ell \in \{2,3\}$.

\begin{pro}\label{pro:Ygroups}
  Let $c \in \N$ with $c \geq 3$, and let $r \in \N_0$.  Then
  $$
  \rk(Y_{c,r}) = 2^{r+1}.
  $$
\end{pro}

\begin{proof}
  We reason similarly as in the proof of
  Proposition~\ref{pro:Xgroups}.  Write $Y := Y_{c,r}$ and argue by
  induction on $r$.  If $r = 0$, the group $Y \cong \SD_{2^{c+1}}$ is
  metacyclic but not cyclic, and the assertion follows.  From now on
  suppose that $r \geq 1$.

  By considering the base subgroup of $Y = \SD_{2^{c+1}} \wr W_r(2)$,
  we obtain
  $$
  \rk(Y) \geq \rk((\SD_{2^{c+1}})^{2^r}) = 2^r \rk(\SD_{2^{c+1}}) =
  2^{r+1},
  $$  
  and it remains to show that $\rk(Y) \leq 2^{r+1}$.  Consider an
  arbitrary subgroup $H \leq Y$, and decompose $Y \cong Y_{c,r-1} \wr
  C_2$ as $Y = B \rtimes C$, where $B = B_1 \times B_2$ with $B_i
  \cong Y_{c,r-1}$, $i \in \{1,2\}$, and $C \cong C_2$.  If $H \leq
  B$, then induction on $r$ yields
  $$
  d(H) \leq \rk(B) = 2 \rk(Y_{c,r-1}) = 2 \cdot 2^r = 2^{r+1}.
  $$
   
  Thus it remains to consider the case where $H$ is not contained in
  $B$.  Decompose $Y \cong \SD_{2^{c+1}} \wr W_r(2)$ as $Y = E \rtimes
  W$, where $E \cong (\SD_{2^{c+1}})^{2^r}$ and $W \cong W_r(2)$.
  Since $\SD_{2^{c+1}}$ is metacyclic, $E$ is metabelian and, in fact,
  decomposes via a characteristic subgroup $F$ into homocyclic parts
  $E/F \cong (C_2)^{2^r}$ and $F \cong (C_{2^c})^{2^r}$.
  
  Choose $h \in H \setminus B$ and regard $(H \cap E)F/F$ and $H \cap
  F$ as $\Z_2 \langle h \rangle$-modules.  In view of the inequality
  $$
  d(H) \leq d(HE/E) + d_{\Z_2 \langle h \rangle}((H \cap E)F/F) +
  d_{\Z_2 \langle h \rangle}(H \cap F),
  $$
  it suffices to establish: (i) $d(HE/E) \leq 2^{r-1}$, (ii)
  $d_{\Z_2 \langle h \rangle}((H \cap E)F/F) \leq 2^{r-1}$ and
  (iii) $d_{\Z_2 \langle h \rangle}(H \cap F) \leq 2^r$.

  Since $HE/E$ embeds into $W \cong X_{1,r-1}(2)$,
  Proposition~\ref{pro:Xgroups} yields~(i).  Considering (ii), we
  write $E = E_1 \times E_2$, where $E_i = E \cap B_i \cong
  (\SD_{2^{c+1}})^{2^{r-1}}$ for $i \in \{1,2\}$.  It is convenient to
  set $K := H \cap E$ and to use $\overline{\cdot}$ to denote
  reduction modulo $F$.  Note that $h$ transposes the two factors in
  the decomposition $\overline{E} = \overline{E_1} \times
  \overline{E_2}$ of the elementary $2$-group $\overline{E}$.  Thus
  $h$ conjugates $\overline{K} \cap \overline{E_2}$ injectively into
  $\overline{E_1} \cong \overline{E}/\overline{E_2}$, and we conclude
  that
  $$
  d_{\Z_2 \langle h \rangle}(\overline{K}) \leq \left(
    d(\overline{K} / \overline{E_2}) - d(\overline{K} \cap
    \overline{E_2}) \right) + d(\overline{K} \cap \overline{E_2}) =
    d(\overline{K} / \overline{E_2}).
  $$
  As $\overline{E} / \overline{E_2} \cong \overline{E_1}$ is an
  elementary $2$-group of rank $2^{r-1}$, this proves~(ii).  Finally,
  $h$ also transposes the two factors $F_i := F \cap E_i$, $i \in
  \{1,2\}$, of $F = F_1 \times F_2$.  Working with the pre-image $M$
  of $H \cap F$ in the base group of the $2$-adic group $(\Z_2 \ltimes
  \Z_2^*) \wr W$, we apply Lemma~\ref{lem:no_fixed_points} to
  deduce~(iii).
\end{proof}

Again, it would be interesting to give a complete description of all
subgroups witnessing the rank of $Y_{c,r}$.  We are ready to prove
Theorem~\ref{thm:GL_d_F_p} and the following complementary result.

\begin{pro}\label{pro:Qp_bound}
  Suppose that $p$ is odd, and let $d \in \N$.  Let $G$ be a maximal
  finite $p$-subgroup of $\GL_d(\Q_p)$.  Then $\rk(G) = \lfloor d
  /(p-1) \rfloor$.
\end{pro}

\begin{proof}[Proof of Theorem~\ref{thm:GL_d_F_p} and
  Proposition~\ref{pro:Qp_bound}]
  The claims are direct consequences of Propositions
  \ref{pro:max_p_subs}, \ref{pro:Xgroups} and~\ref{pro:Ygroups},
  combined with the description of Sylow subgroups of finite symmetric
  groups in terms of wreath products.
\end{proof}


\section{Bounding the $\ell$-rank}

\begin{pro}\label{pro:key}
  Suppose that $p$ is odd.  Let $G$ be a finite $p$-group, and let $M$
  be a faithful $\Z_p G$-module which is free and of finite rank as a
  $\Z_p$-module.  Then
  $$
  d(G) + d_{\Z_p G}(M) \leq \dim_{\Z_p}(M).
  $$
\end{pro}

\begin{proof}
  Clearly, we may suppose that $G$ is not trivial.  Let $C = \langle x
  \rangle$ be a central subgroup of $G$ of order $p$.  According to
  \cite[Theorem~2.6]{HeRe62}, there are three types of indecomposable
  $\Z_p C$-modules: they are the trivial module $\Z_p$ of dimension
  $1$, the free module $\Z_p C$ of dimension~$p$, and the
  $(p-1)$-dimensional module $\Z_p C/\Phi_p(x) \Z_p C$ where $\Phi_p =
  1 + X + \ldots + X^{p-1}$.  As a $\Z_p C$-module, $M$ is the direct
  sum of indecomposable modules of these three types.

  Set $N := \Cen_M(C) = \{ m \in M \mid m^x = m \}$.  Since $C$ is
  central in $G$ and since $C$ acts faithfully on $M$, the set $N$
  constitutes a proper $\Z_p G$-submodule of $M$.  We distinguish two
  cases: $N = \{ 0 \}$ and $N \not = \{ 0 \}$.

  First suppose that $N \not = \{ 0 \}$.  In this case, $G$ acts
  naturally on the quotient module $M_1 := M/N$ with image $G_1 \leq
  \GL(M/N)$ and kernel $G_2$, say.  The finite group $G_2$ acts
  faithfully on $M_2 := N$; indeed, the $\Q_p G$-module $\Q_p
  \otimes_{\Z_p} M$ is completely reducible and $\Q_p \otimes_{\Z_p}
  M_2$ admits a complement isomorphic to $\Q_p \otimes_{\Z_p} M_1$.
  Clearly, $d(G) \leq d(G_1) + d(G_2)$ and $d_{\Z_p G}(M) \leq d_{\Z_p
    G_1}(M_1) + d_{\Z_p G_2}(M_2)$.  Therefore induction gives
  \begin{align*}
    d(G) + d_{\Z_p G}(M) & \leq d(G_1) + d(G_2) + d_{\Z_p G_1}(M_1) +
    d_{\Z_p G_2}(M_2) \\
    & \leq  \dim_{\Z_p}(M_1) + \dim_{\Z_p}(M_2) \\
    & = \dim_{\Z_p}(M).
  \end{align*}

  Now suppose that $N = \{0\}$.  Then, as a $\Z_p C$-module, $M$
  necessarily decomposes into a direct sum of indecomposable $\Z_p
  C$-modules each of which is isomorphic to $\Z_p C/\Phi_p(x) \Z_p C$.
  Clearly,
  $$
  d_{\Z_p G}(M) \leq d_{\Z_p C}(M) = (p-1)^{-1} \dim_{\Z_p}(M).
  $$
  On the other hand, Proposition~\ref{pro:Qp_bound} shows that $d(G)
  \leq (p-1)^{-1} \dim_{\Z_p}(M)$ and thus
  $$
  d(G) + d_{\Z_p G}(M) \leq 2(p-1)^{-1} \dim_{\Z_p}(M) \leq
  \dim_{\Z_p}(M),
  $$
  as wanted.
\end{proof}

\begin{lem}\label{lem:witnesses}
  Let $G$ be an insoluble pro-$p$ group, and suppose that $p$ is odd.
  Then for every open normal subgroup $N$ of $G$ there exists an open
  normal subgroup $K$ of $G$ such that $K \subsetneqq N$ and $d(K)
  \geq d(N)$.
\end{lem}

\begin{proof}
  Clearly, the assertion holds if $G$ is not finitely generated.  Now
  suppose that $d(G) < \infty$, and let $N$ be an open normal subgroup
  of~$G$.  Let $\Phi(N) = N^p[N,N]$ denote the Frattini subgroup of
  $N$.  Since $G$ is insoluble, we have $M := [\Phi(N),N] \neq 1$.
  Observe that $M$ is normal in $G$, but not necessarily open.
  Nevertheless, working modulo an open normal subgroup $U$ of $G$ such
  that $U \subsetneqq MU \subseteq \Phi(N)$, we find an open normal
  subgroup $L$ of $G$ such that $L \subseteq MU \subseteq \Phi(N)$ and
  $MU/L \cong C_p$.  Now an argument similar to the one given in the
  proof of \cite[Theorem~3.3]{GoKl10} shows that
  $$
  K := \{ g \in N \mid [\Phi(N),g] \subseteq L \}
  $$
  satisfies $d(K) \geq d(N)$.  One checks easily that $K
  \trianglelefteq G$ and $K \subsetneqq N$.
\end{proof}

Recall that $\pi(G)$ denotes the periodic radical, viz.\ the unique
maximal finite normal subgroup, of a profinite group $G$ of finite
rank.  

\begin{lem}\label{lem:pi}
  Let $G$ be a profinite group of finite rank and $H$ a subgroup of
  $G$.  If $H$ is normal or open in $G$, then $\pi(H) \leq \pi(G)$.
\end{lem}

\begin{proof}
  In both cases, the union of all conjugates of $\pi(H)$ in $G$ is a
  finite normal subset of $G$ consisting of elements of finite order.
  By Dicman's Lemma, this set generates a finite normal subgroup of
  $G$.
\end{proof}

For the next two results, we denote by $d_\ell(G)$ the minimal
cardinality of a generating set for a Sylow pro-$\ell$ subgroup of a
profinite group $G$.

\begin{pro}\label{pro:bound_for_dp}
  Let $G$ be a compact $p$-adic Lie group such that $\pi(G) = 1$, and
  suppose that $p$ is odd.  Then $d_p(G) \leq \dim(G)$.
\end{pro}

\begin{proof}
  By Lemma~\ref{lem:pi}, we may assume without loss of generality that
  $G$ is a pro-$p$ group.  If $\dim(G) = 0$, then $G = \pi(G) = 1$ and
  there is nothing further to prove.  Hence suppose that $\dim(G) \geq
  1$.  Choose a normal subgroup $N$ of $G$ such that $G/N$ is just
  infinite.  Note that $\pi(G/N) = 1$ and $\pi(N) = 1$, by
  Lemma~\ref{lem:pi}.  Thus, if $N \neq 1$ we apply induction to
  deduce that
  $$
  d(G) \leq d(G/N) + d(N) \leq \dim(G/N) + \dim(N) = \dim(G).
  $$
  It remains to consider the case that $N = 1$, i.e.\ that $G$ is just
  infinite.

  First suppose that $G$ is soluble.  Then $G$ is virtually abelian;
  we find an abelian open normal subgroup $A$ of $G$.  Put $C :=
  \Cen_G(A)$.  Then $\lvert C : \Zen(C) \rvert \leq \lvert C : A
  \rvert < \infty$, and hence $[C,C]$ is finite.  Since $G$ is just
  infinite, we conclude that $C$ is abelian and self-centralising.
  Thus Proposition~\ref{pro:key} shows that $d(G) \leq \dim(C) =
  \dim(G)$.

  Next suppose that $G$ is insoluble.  For a contradiction assume that
  $d(G) > \dim(G)$.  By Lemma~\ref{lem:witnesses}, we find a strictly
  descending chain
  $$
  G = N_1 \supset N_2 \supset N_3 \supset \ldots
  $$
  of open normal subgroups $N_i$ of $G$ such that $d(N_i) \geq d(G)$
  for all $i \in \N$.  Fix a uniformly powerful open normal subgroup
  $U$ of $G$.  Since $\rk(U) = \dim(U) < d(G)$, we have $N_i \not
  \subseteq U$ for each $i \in \N$, witnessed by $g_i \in N_i
  \setminus U$, say.  Since $G$ is compact, there is a subsequence of
  $g_i$, $i \in \N$, which converges to an element $g \in G \setminus
  U$.  Clearly, we also have $g \in \bigcap \{ N_i \mid i \in \N \} =:
  M$.  Since $M$ is a normal subgroup of infinite index in $G$ and
  since $G$ is just infinite, we conclude $g \in M = \{ 1 \} \subseteq
  U$, a contradiction.
\end{proof}

\begin{pro} \label{pro:bound_for_dl}
  Let $G$ be a compact $p$-adic Lie group such that $\pi(G) = 1$, and
  let $d := \dim(G)$.  Then for every prime $\ell$ with $\ell \not =
  p$,
  $$
  d_\ell(G) \leq \rk_\ell(\GL_d(\Z_p)) = \rk_\ell(\GL_d(\F_p)).
  $$
\end{pro}

\begin{proof}
  Arguing as at the beginning of the proof of
  Proposition~\ref{pro:bound_for_dp}, we may suppose that $G$ is just
  infinite.  Let $S$ be a Sylow pro-$\ell$ subgroup of~$G$, and choose
  a uniformly powerful open normal pro-$p$ subgroup $U$ of~$G$.  Then
  $S$ is a finite $\ell$-group, and we claim that $S$ acts faithfully
  on $U$ by conjugation.  Indeed, since $G$ is just infinite,
  $\Cen_G(U)$ is either trivial or open in $G$.  In the former case
  there is nothing further to show.  In the latter case, writing $C :=
  \Cen_G(U)$, we note that $\lvert C : \Zen(C) \rvert \leq \lvert C :
  C \cap U \rvert < \infty$ so that $[C,C]$ is a finite normal
  subgroup of $G$.  Since $G$ is just infinite, we conclude that $C$
  is abelian and pro-$p$.  Hence $S \cap C = 1$, and again $S$ acts
  faithfully on $U$.

  Since $U$ is uniformly powerful, the set $U$ admits the structure of
  a $\Z_p$-Lie lattice with additive group isomorphic to $\Z_p^d$,
  where $d := \dim(U) = \dim(G)$, and the faithful action of $S$ on
  $U$ translates into an embedding of $S$ into $\GL_d(\Z_p)$.
  Dividing out the first principal congruence subgroup
  $\GL_d^1(\Z_p)$, we obtain an embedding of $S$ into $\GL_d(\F_p)$.
\end{proof}

\begin{proof}[Proof of Theorem~\ref{thm:A}]
  Let $G$ be a compact $p$-adic Lie group such that $\pi(G) = 1$, and
  suppose that $p$ is odd.  The assertion about $\rk_2(G)$ follows
  directly from Proposition~\ref{pro:bound_for_dl} and
  Theorem~\ref{thm:GL_d_F_p}.  For the assertions about $\rk_\ell(G)$,
  $\ell > 2$, it suffices to show that
  \begin{enumerate}
  \item[(i)] $\rk_p(G) \geq \dim(G)$,
  \item[(ii)] $\rk_\ell(G) \leq \dim(G)$ for every odd prime $\ell$.
  \end{enumerate}
  For every uniformly powerful open pro-$p$ subgroup $U$ of $G$ one
  has $d(U) = \dim(U) = \dim(G)$, and this implies (i).  Let $\ell$ be
  an odd prime.  It is known that
  $$
  \rk_\ell(G) = \sup \{ d_\ell(H) \mid H \leq G \text{ an open
    subgroup} \};
  $$
  cf.\ \cite[Proposition 3.11]{DiDuMaSe99}.  Thus (ii) follows from
  Lemma~\ref{lem:pi} and
 Propositions~\ref{pro:bound_for_dp} and
  \ref{pro:bound_for_dl}, in conjunction with
  Theorem~\ref{thm:GL_d_F_p}.
\end{proof}


\section{Bounding the rank}

We will make use of the following theorem which was proved for soluble
groups by Kov\'acs, and for all finite groups (using the
classification of finite simple groups) by Lucchini~\cite{Lu89} and
Guralnick~\cite{Gu89}.

\begin{thm}
  If every Sylow subgroup of a finite group $G$ can be generated by
  $d$ elements, then $d(G) \leq d+1$.
\end{thm}

\begin{cor} \label{cor:kovacs}
 For every profinite group $G$ one has 
 \begin{equation}\label{equ:lucchini}
   \rk(G) \leq \sup \{ \rk_\ell(G) \mid \text{$\ell$ prime} \} + 1.
 \end{equation}
\end{cor}

\begin{proof}[Proof of Theorem~\ref{thm:B}]
  Claims (1) and (3) are direct consequences of
  Corollary~\ref{cor:kovacs} and Theorem~\ref{thm:A}.  To deduce
  claims (2) and (4) we use the additional analysis given
  in~\cite{Lu90}, which distinguishes between groups of zero and
  non-zero presentation rank.

  First suppose that $G$ is prosoluble.  Then all its finite quotients
  have zero presentation rank.  Hence \cite[Theorem~1]{Lu90} shows
  that equality in \eqref{equ:lucchini} implies $\rk(G) \leq
  \rk_\ell(G) + 1$ for some odd prime $\ell$.  Thus claim (2) follows
  from Corollary~\ref{cor:kovacs} and Theorem~\ref{thm:A}.

  Finally suppose that $G$ has trivial $(p-1)$-torsion and without
  loss of generality assume that $\dim(G) \geq 1$.  Then $\rk_p(G) =
  \dim(G)$ and $\rk_\ell(G) = 0 \leq \dim(G) -1$ for primes $\ell$
  with $p \equiv_\ell 1$.  Furthermore,
  Proposition~\ref{pro:bound_for_dl} and Theorem~\ref{thm:GL_d_F_p}
  show that $\rk_\ell(G) \leq \lfloor \dim(G)/2 \rfloor \leq \dim(G) -
  1$ for all remaining primes $\ell$.  Corollary~\ref{cor:kovacs}
  gives $\rk(G) \in \{ \dim(G), \dim(G)+1 \}$, and \cite[Theorems~1
  and~2]{Lu90} rule out the possibility $\rk(G) = \dim(G)+1$.  This
  proves claim (4).
\end{proof}


\bigskip

\noindent
\textbf{Acknowledgements.} The author thanks Frank Calegari for posing
a question which led to the present article.  He is also grateful to
Gabriele Nebe for pointing out a gap in the arguments given in an
earlier version of this paper.


\end{document}